\documentstyle[12pt]{article}
\newcommand{\nt}{\noindent}

\begin{document}
\begin{center}
{\bf  An Application of a Log Version of the Kodaira Vanishing Theorem 
to Embedded Projective Varieties}

\medskip

by

\medskip

Aaron Bertram\footnote
{Partially supported by a Sloan research fellowship.}

\end{center}

\bigskip

\nt {\bf 0. Introduction.} 
Let $Y \subset {\bf CP}^n$ be a 
smooth complex projective subvariety of codimension $r$, and 
let ${\cal I}_Y$ be the ideal sheaf of the embedding, with 
${\cal I}_Y^k \subset {\cal O}_{{\bf P}^n}$ denoting its $k$th
power. In this paper, we will be ineterested in the following two
integer invariants of the embedding:

\medskip

$d_Y$, the minimum of the degrees $d$ such that $Y$ 
is a scheme-theoretic 
intersection of hypersurfaces of degree at most $d$, 
and given $d_Y$,

\medskip

$e_Y$, the minimum of the integers $e$ such that:
$$H^i({\bf P}^n,{\cal I}_Y^k(p)) = 0 \ \mbox{for all $i > 0,k \ge 0$
and $p \ge e + (k-1)d_Y$}.$$

An upper bound for $e_Y$ (which is sharp if $Y$ is a complete intersection)
was computed by the author in collaboration with Ein and Lazarsfeld:

\medskip

\nt {\bf Theorem 1:}(\cite{BEL}, Proposition 1) Suppose $Y$ is 
scheme-theoretically cut out by equations of degrees
$$d_Y = d_1 \ge d_2 \ge ... \ge d_m.$$
Then $e_Y \le d_1 + ... + d_r - n$. 
(Recall that $r$ is the codimension of $Y$.)

\medskip

The idea in this paper is to show how a generalized
``log'' version of the Kodaira vanishing theorem
can be employed to improve the results of 
Theorem 1 (which was also proved by
Kodaira vanishing) when we have more
knowledge about the equations for $Y$.
The idea is to
find a hypersurface $F \subset {\bf P}^n$ which has 
high multiplicity along $Y$, is ``log canonical'' near $Y$, and
has relatively small degree, then to invoke Kodaira vanishing
on the blow-up of ${\bf P}^n$ along $Y$. In the context of 
Theorem 1, the hypersurface $F$ is approximately a divisor
with normal crosssings (see its proof in \S 2).
However one of the main points of this paper is the observation
that even in the most familiar of projective embeddings, 
log canonical divisors quite different from normal crossings
divisors seem to play an important role.  
  
\medskip

All the new cases we consider are determinantal, 
in the sense that $d_Y = 2$ and a collection of quadrics which 
scheme-theoretically cut out 
$Y$ arise either as $2\times 2$ minors or $4\times 4$ Pfaffians of 
a matrix of linear forms on ${\bf P}^n$. In each of these cases,
the hypersurface $F$ is constructed out of minors or 
Pfaffians of all sizes.  
The results are most satisfactory for the``universal'' determinantal varieties
where the key is to observe that the theories of complete 
linear maps and quadrics (and a version involving Pfaffians in
the skew  case, which seems not to have been 
previously worked out) give us the information we need to check
whether hypersurfaces built out of minors have 
mild enough singularities. 
We apply the same idea to curves embedded by a line bundle of 
large degree, obtaining similarly satisfactory results in genus $0$ and $1$.
However in higher genus, some complications arise, and the 
results obtained here are probably not the best.

The contents of the paper are as follows. In \S 1, we review some of the 
relevant definitions and results of the log minimal 
model program leading up
to the log version of the Kodaira vanishing
theorem, due to Nadel. In addition, we state a useful Bertini
property, due to Koll\'ar.
In \S 2 we explain how a hypersurface
$F \subset {\bf P}^n$  
yields an upper bound
on $e_Y$, and use it in subsequent sections to give:

\bigskip

(i) A reproof of Theorem 1 taking
$F$ to be a sum of hypersurfaces 
of degree $d_i$ (for $1 \le i \le r$) 
which are general among those vanishing on $Y$.

\bigskip

(ii) An upper bound for 
$e_Y$ which is independent of the dimensions of vector spaces $V$ and $W$
for each of the three universal determinantal varieties:

\medskip

\hskip .2in (a) (Generic) $Y = {\bf P}(V) \times {\bf P}(W)$, Segre embedding.
$e_Y \le -1$.

\medskip

\hskip .2in (b) (Symmetric) $Y = {\bf P}(V)$, quadratic Veronese embedding. 
$e_Y \le 0$.  

\medskip

\hskip .2in (c) (Skew) $Y = G(2,V)$, Pl\"ucker embedding. $e_Y \le -3$.

\bigskip

(iii) A proof that $e_Y \le 1$ when $Y = C$ is a Riemann surface
of genus $g$ embedded by a complete linear series in the following cases:

\medskip

\hskip .2in (a) $g = 0$ or $1$.

\medskip

\hskip .2in (b) the degree of the embedding is sufficiently large 
($> \frac {8g + 2}3$ will do).

\medskip

\hskip .2in (c) the degree of the embedding is at least 
$2g+3$, but with ``gaps''.

\medskip

\nt (Notice that in (ii) and (iii), Theorem 1 would give only 
$e_Y \le 2r - n$.)

\bigskip

\nt {\bf Remarks:} 
In \cite{W}, Wahl
proves the vanishing H$^1({\bf P}^n,{\cal I}_X^2(p)) = 0$ for all
$p\ge 3$ and $X \subset {\bf P}^n$ embedded by a complete linear
series in the following cases: 

\medskip

(1) $X$ is projective space,

\medskip

(2) $X$ is arbitrary, but the linear series is sufficiently
ample, and

\medskip

(3) $X$ is a general canonical curve of genus $\ge 3$.

\medskip

Since in all these cases the embedding is projectively normal
and $X$ is scheme-theoretically cut out by quadrics, Wahl's
results may be a special case of the more general 
property $e_Y \le 1$. This we know to be the case for (1) and
(2)  when $X$ is a curve by the results of \S 4. It would
be interesting to know  whether or not this property does indeed
hold in this generality (as well as the case of an embedding of a
curve of degree $2g+3$ or more). 

\medskip

The next remark is more of a confession, really. The
invariant $e_Y$ defined here probably ought to be 
modified to conform with Castelnuovo-Mumford regularity. Recall
that  if ${\cal F}$ is a sheaf on ${\bf P}^n$, then ${\cal F}$
is defined to be $m$-regular if H$^i({\bf P}^n,{\cal F}(m-i)) =
0$ for all $i > 0$. The main feature of regularity is that
$m$-regular implies $m+1$-regular. But of course this pattern
of vanishing in case ${\cal F}$ is a power of 
the ideal sheaf of $X$ does not lend
itself to proof by vanishing theorems as outlined in this paper.

\medskip 
 
On the other hand, in \cite{T} \S 6, Thaddeus obtains some
similar vanishing results using ordinary Kodaira
vanishing in case $X$ is a
curve embedded by a line bundle of large degree. In his case, 
the vanishing takes place on spaces obtained from 
$Y$ by a sequence of flips. These flips are shown to 
be log flips in the sense of the log minimal model program
in \cite{B2} using log canonical divisors of precisely the sort 
we use here to prove vanishing. Probably the sharpest
results would be obtained by applying the generalized
Kodaira vanishing theorem on these flipped spaces and
transferring the vanishing results back to $Y$. It seems entirely
possible that such a procedure will yield a pattern of 
vanishing which does conform with Castelnuovo-Mumford
regularity. 
 
\newpage

\nt {\bf \S 1. Log Kodaira Vanishing:} 
Let $X$ be a smooth complex projective variety of dimension $n$. 
The following definitions 
are standard to the experts,  
but are perhaps not widely known:

\bigskip

\nt {\bf Definitions:} (a) A finite ${\bf Q}$-linear combination
$F = \sum \alpha_iF_i$ of distinct prime divisors of $X$ is
called a 
${\bf Q}$-divisor. It is {\bf effective} 
if each $\alpha_i \ge 0$. 
Intersection with divisors extends by linearity to give 
well-defined rational numbers $F^n$ and $F.B$,
given a ${\bf Q}$-divisor $F$ and a curve $B \subset X$. 
In particular, numerical equivalence
extends to an equivalence relation on ${\bf Q}$-divisors.

\medskip

(b) A(n equivalence class of) ${\bf Q}$-divisor(s) $A$ is {\bf nef and big} if:

\medskip

\hskip .5in (i) $A.B \ge 0$ for all curves $B \subset X$ and

\medskip

\hskip .5in (ii) $A^n > 0$.

\bigskip

Given an effective ${\bf Q}$-divisor 
$F = \sum \alpha_iF_i$ and a birational
morphism $f:\widetilde X \rightarrow X$, 
let $E$ be the $f$-exceptional
divisor, and let $\{E_j\}$ be the components of $E$. 
Also let $f^*(F_i)$ and $f_*^{-1}(F_i)$ be the total and strict transforms,
respectively, of $F_i$
on $\widetilde X$, extending the usual notions by
linearity. 

\medskip

(c) If $f:\widetilde X \rightarrow X$ has the property that
$\widetilde X$ is smooth and 
$\sum E_j + \sum f_*^{-1}(F_i)$
is a normal crossings divisor with smooth components, then  
$f$ is called a {\bf log resolution} of the pair $(X,F)$.
Given such an $f$, one attaches a rational number to each $E_j$ 
and $f_*^{-1}(F_i)$, called the
{\bf discrepancy} of $f$, as follows:

\medskip

\hskip .2in (i) The discrepancy at $f_*^{-1}(F_i)$ is $-\alpha_i$. 

\medskip

\hskip .2in (ii) The discrepancy at $E_j$ is its coefficient in the difference:
$$(K_{\widetilde X} + f_*^{-1}(F)) - f^*(K_X +F)$$

(d) Given an effective ${\bf Q}$-divisor $F$ and a log resolution
$f:\widetilde X \rightarrow X$,

\bigskip

{\bf discrep}$(X,F,f)$ is the minimum of the discrepancies of type (ii), and

\medskip

{\bf totaldiscrep}$(X,F,f)$ is the minimum of all the discrepancies. 

\bigskip

\nt {\bf Remark:} We have limited ourselves here to 
smooth $X$, since that is all we will need to consider. 
See \cite{Ketal} or \cite{K} for the general definitions 
when $X$ is not assumed to be smooth, as well as a proof of the following: 

\bigskip

\nt {\bf Basic Observation:} The following definitions are intrinsic to a 
pair $(X,F)$ (i.e. they do not depend upon the log resolution $f:\widetilde X \rightarrow X$):

\bigskip

\nt $(X,F)$ is {\bf log canonical} (or lc) if
totaldiscrep$(X,F,f)
\ge -1$.

\bigskip

\nt $(X,F)$ is {\bf Kawamata log terminal} (or klt) if
totaldiscrep$(X,F,f) > -1$.

\bigskip

\nt $(X,F)$ is {\bf purely log terminal} (or plt) if it is
log canonical, and if, in addition,
discrep$(X,F,f) > -1$.

\bigskip

\nt {\bf Examples:} (1) $(X,F)$ is log canonical when $F$ is a Cartier divisor
with smooth components and normal crossings. (the identity is a log resolution
of $(X,F)$,
and all the discrepancies are $-1$ or $0$.)

\medskip

(2) Suppose $Z_1 \subset Z_2 \subset ... \subset Z_k \subset X$
are closed subvarieties, 
and $F$ is an effective (Cartier) divisor 
on $X$. Suppose blowing up the strict transforms of each 
$Z_j$ in order is 
a sequence of blow-ups along smooth centers so that the composition
of blow-downs $f:\widetilde X \rightarrow X$ is a log resolution of $(X,F)$.
Let $m_j$ be the multiplicity of $F$ at the generic point of $Z_j$.
Then:

\medskip

\hskip .5in $(X,F)$ is lc if $m_j \le \mbox{codim}_X(Z_j)$ for all $j$ and

\medskip

\hskip .5in $(X,F)$ is plt if $m_j < \mbox{codim}_X(Z_j)$ 
(it is only klt if $F = \emptyset$(!)). 

\medskip

(This is an immediate consequence of Riemann-Hurwitz.)

\medskip

(3) Given a log canonical pair $(X,F)$ and a rational number  
$0 < \epsilon < 1$, then
the pair $(X,(1-\epsilon)F)$ is klt. (Immediate from the definitions.)

\medskip

\nt {\bf More Definitions:} Let $F$ be an effective
${\bf Q}$-divisor on $X$. For each $x\in X$, one says $(X,F)$ is {\bf not 
lc} (resp. {\bf not klt}) {\bf at $x$} if there is a subvariety 
$x\in Z \subset X$, a log resolution $f:\widetilde X \rightarrow X$ and 
an exceptional (or strict-transform) divisor $E_Z \subset \widetilde X$ such that
$f(E_Z) = Z$ and the discrepancy at $E_Z$ is $< -1$ (resp. $\le -1$). The
following subsets of $X$ are known to be closed:

\medskip

{\bf Nklt}$(X,F) := \{\ x\in X \ | \ (X,F) \ 
\mbox{is not klt at $x$}\}$, and 

\medskip

{\bf Nlc}$(X,F) := \{\ x\in X \ | \ (X,F) \ \mbox{is not lc at $x$}\}$
(closed by (3) above).

\bigskip

We will use the following very simple case of a Bertini property 
due to Koll\'ar which tells us that the Nlc and Nklt loci for general
members of linear series 
can be detected ``pointwise'' (again, see \cite{K} for a much more
general version).   

\medskip

Suppose $F$ is an effective ${\bf Q}$-divisor and 
$|B_1|,...,|B_k|$ are linear series on $X$.  Let 
$B_i^g$ denote a general member of $|B_i|$, let $B^g := B_1^g + ... + B_k^g$,
and let $b_1,...,b_k$ be rational numbers between $0$ and $1$. Then using
(4.8.1-2) of \cite{K}, we obtain:

\medskip

\nt {\bf Bertini Property:} If $x\not\in \ \mbox{Nlc}(X,F+B^g)$ for each 
$x$ in some subset $W\subset X$, then 
Nlc$(X,F+B^g) \cap W = \emptyset$. The 
same is true with lc replaced by klt provided that the $b_i$ 
are strictly less than $1$.

\bigskip

(The point is that a priori the choice of $B^g$ 
could depend upon $x$.)

\bigskip

\nt {\bf Example:} If $(X,F)$ is log canonical and the $|B_i|$ are
all base-point-free, then the Bertini property shows that
$(X,F+\sum_{i=1}^k B_i^g)$ is log canonical.

\bigskip

The following theorem is due to Alan Nadel 
(see Koll\'ar's notes, Theorem 2.16 for a more general version when
$X$ is  allowed some singularities).

\bigskip

\nt {\bf Theorem (Log Kodaira Vanishing):} Suppose that   
$F$ is an effective ${\bf Q}$-divisor on $X$, 
$A$ is another
${\bf Q}$-divisor which is nef and big, and that
$L$ is a line bundle on $X$ satisfying:
$$L \equiv K_X + F + A.$$

Then there is an ideal sheaf ${\cal J}$ on X (called Nadel's multiplier 
ideal sheaf) with the following properties:

\medskip

(i) ${\cal O}_X/{\cal J}$ is supported 
on Nklt$(X,F)$ (which is therefore closed!), and 

\medskip (ii) $H^i(X,{\cal J} \otimes L) = 0 \ \mbox{for all} \ i > 0.$ 

\bigskip

And the obvious corollary:

\bigskip

\nt {\bf Corollary:} If $(X,F)$ is klt in the theorem, then:
$$H^i(X, L) = 0 \ \mbox{for all} \ i > 0.$$

\newpage

\nt {\bf \S 2. The Strategy (and Reproof of Theorem 1):} 
We return now to the set-up from the introduction. $Y \subset {\bf P}^n$ is 
a smooth projective subvariety of codimension $r > 0$. Let:
$$X := \ \mbox{bl}\ ({\bf P}^n,Y), \ \mbox{the blow-up of ${\bf P}^n$ along $Y$,}$$
and let $H$ and $E$ be hyperplane 
and exceptional divisors on $X$

\medskip

Here are a few standard observations about $X$:

\medskip

\hskip .2in (1) $K_X \equiv -(n+1)H + (r-1)E$ (Riemann-Hurwitz).

\medskip

\hskip .2in (2) $H - \epsilon E$ is ample for $0 < \epsilon << 1$ (Kleiman's
Criterion).

\bigskip

Our strategy for seeking upper bounds for $e_Y$ rests on the 
following proposition, which is the essential observation of the 
paper.

\bigskip

\nt {\bf Proposition 2.1:} If there is an effective ${\bf Q}$-divisor $F$ on $X$ such that:

\medskip

(i) $F \equiv (e+n)H - rE$ and 

\medskip

(ii) Nlc$(X,F) \cap E = \emptyset$,

\medskip

\nt then $e_Y \le e$.

\medskip

{\bf Proof:} Given such an $F$, then for each $\epsilon \in {\bf Q}$
satisfying  $0 < \epsilon < 1$, we
would have Nklt$(X,(1-\epsilon)F) \cap E = \emptyset$, and using (1), 
$$pH - E \equiv K_X + (1-\epsilon)F + A$$
where $A \equiv (p + 1 - e + (e+n)\epsilon)H - r\epsilon E$. 
If additionally, $\epsilon << 1$, then by (2), $A$ is ample (hence nef and big)
provided that $p \ge e$.

\medskip

Moreover, $d_YH - E$ (and all positive multiples) is base-point-free on $X$, 
by definition of $d_Y$, so that
for each positive integer $k$,
$$pH - kE \equiv K_X + (1-\epsilon)F + A$$  
where $A \equiv (p+1 - e-(k-1)d_Y + (e+n)\epsilon)H - r\epsilon E + 
(k-1)(d_YH-E)$ is ample provided that 
$p \ge e+(k-1)d_Y$. 

Thus, the log Kodaira vanishing theorem (using (ii)) tells us that
$$H^i(X,{\cal J}\otimes {\cal O}_X(pH-kE)) = 0 \ \mbox{for all} \ i > 0, p \ge
e + (k-1)d_Y $$
where ${\cal J}$ is an ideal sheaf on $X$ with the property that 
the support of ${\cal O}_X/{\cal J}$ is disjoint from $E$.

\medskip

We can therefore identify ${\cal J}$ with 
its direct image in ${\bf P}^n$, and it is a consequence of 
the theorem of formal functions (\cite{H},III.11) that:
$$H^i({\bf P}^n,{\cal J} \otimes {\cal I}_Y^k \otimes {\cal O}_{{\bf P}^n}(p))
= 0 \ \mbox{for all} \ i > 0, p \ge
e + (k-1)d_Y.$$

To conclude the vanishing without ${\cal J}$, we use the 
fact that the ideal sheaves ${\cal J}$ and ${\cal I}^k_Y$ have disjoint
cosupport to conclude that ${\cal O}_{{\bf P}^n}/{\cal J}$ is a direct summand 
of ${\cal O}_{{\bf P}^n}/{\cal J}{\cal I}_Y^k = 
{\cal O}_{{\bf P}^n}/{\cal J}\otimes {\cal I}_Y^k$.  We also use the 
disjoint cosupport in the first of the following 
two exact sequences:
$$0 \rightarrow {\cal J}\otimes {\cal I}_Y^k\otimes {\cal O}_{{\bf P}^n}(p) 
\rightarrow {\cal I}_Y^k \otimes {\cal O}_{{\bf P}^n}(p) 
\rightarrow {\cal O}_{{\bf P}^n}(p)/{\cal J}
\rightarrow 0,$$
$$0 \rightarrow {\cal J}\otimes {\cal I}_Y^k\otimes {\cal O}_{{\bf P}^n}(p) 
\rightarrow {\cal O}_{{\bf P}^n}(p) 
\rightarrow {\cal O}_{{\bf P}^n}(p)/({\cal J}{\cal I}_Y^k)
\rightarrow 0.$$

If it were the 
case that $H^i({\bf P}^n,{\cal I}_Y^k\otimes {\cal O}_{{\bf P}^n}(p)) \ne 0$,
then from the long exact sequence on cohomology
associated to these two short exact sequences and the 
vanishing above,
we would have 
$H^i({\bf P}^n, {\cal O}_{{\bf P}^n}(p)/{\cal J}) \ne 0$,
hence  
$H^i({\bf P}^n,{\cal O}_{{\bf P}^n}(p)/({\cal J}{\cal I}_Y^k)) \ne
0$ and 
$H^i({\bf P}^n,{\cal O}_{{\bf P}^n}(p)) \ne 0$, a contradiction.

\bigskip

With this proposition, we now have a very fast

\bigskip

{\bf Proof of Theorem 1:} Let $I \subset {\bf C}[x_0,...,x_n]$ be the 
ideal generated by 
the given homogeneous polynomials of degree 
$d_Y = d_1,...,d_m$
which scheme-theoretically cut out $Y$. 
For each $i=1,...,r$, let $I_{d_i}$ be the homogeneous part
of degree $d_i$, and let $|B_i|$ be the corresponding
sub-linear  series of $|d_iH - E|$ on $X$. It follows that 
for each $x\in E$, the sum of general elements  
$B^g := B_{1}^g+ ... + B_{r}^g$ is a normal-crossings 
divisor, hence log canonical at $x$. Thus the Bertini property tells us  
that Nlc$(X,B^g)
\cap E = \emptyset$, so $F:= B^g$ satisfies the
conditions of Proposition 2.1 with 
$e = d_1 + ... + d_r - n$.

\bigskip

\nt {\bf Remark:} 
This proof is essentially the same as the proof in
\cite{BEL}. However, by
making the dependence upon a suitable hypersurface $F$
explicit in Proposition 2.1, a general strategy has
emerged which was not apparent in
\cite{BEL}. Namely, given an embedding $Y \subset {\bf
P}^n$, one wants to find hypersurfaces
which  are highly singular along $Y$ relative 
to their degree but whose strict transform on 
$X$ is log canonical near $E$ . We will see
in the next sections that certain determinantal 
varieties fit nicely into this strategy. 
  
\bigskip

\nt {\bf \S 3. Universal Determinantal Varieties:} Let $V$ and $W$
be vector spaces of dimension $k$ and $m$ respectively, suppose
that $k \le m$, and let $Y$ be one of the following:

\medskip

(a) ${\bf P}(V)\times {\bf P}(W)$ embedded in ${\bf P}^n
:= {\bf P}(V\otimes W)$ ($n = km-1$)
by the Segre embedding, 

\medskip

(b) ${\bf P}(V)$ embedded in ${\bf P}^n := {\bf P}(\mbox{Sym}^2(V))$ 
($n = \frac 12(k^2 + k) - 1$) by the 
quadratic Veronese embedding, or

\medskip

(c) $G(2,V^*)$ embedded in ${\bf P}^n :=  {\bf P}(\wedge^2(V))$ ($n = 
\frac 12(k^2 - k) - 1$) by the Pl\"ucker embedding.

\medskip

Then $Y$ is the rank one locus of a universal
map $\phi$ of vector bundles.
$\phi: V\otimes {\cal O}_{{\bf P}^n} \rightarrow
W^* \otimes {\cal O}_{{\bf P}^n}(1)$ in case (a), and   
$\phi:V\otimes {\cal O}_{{\bf P}^n} \rightarrow
V^* \otimes {\cal O}_{{\bf P}^n}(1)$ in cases
(b) and (c), where $\phi$ is, respectively, symmetric and
skew symmetric. 
Alternatively, one can, of course, choose bases for $W$ and $V$ and 
think of $\phi$ as 
a matrix of linear forms. In each case, $Y$ is the last of a nested
sequence of degeneracy loci in ${\bf P}^n$ determined by
the map
$\phi$. In (a) and (b), let $\Delta_i \subset {\bf P}^n$ be the zero locus 
of  $\wedge^i \phi$, while in (c), let $\Delta_i$ be the zero locus
of $\wedge^{2i}\phi$, (to get this right 
scheme-theoretically, one needs to take the ``square root'' of this map...see
below). Then it is a standard fact that:

\medskip

(a) $Y = \Delta_2 \subset \Delta_3 \subset ... \subset \Delta_k$

\medskip

\nt and each $\Delta_i$ is irreducible of codimension $(k-i+1)(m-i+1)$ 
in ${\bf P}^n$,

\medskip

(b) $Y = \Delta_2 \subset \Delta_3 \subset ... \subset \Delta_k$

\medskip

\nt and each $\Delta_i$ is irreducible of codimension $\left( {k-i+2}
\atop 2\right)$ in ${\bf P}^n$,

\medskip

(c) $Y = \Delta_2 \subset \Delta_3 \subset ... \subset \Delta_{[\frac k2]}$

\medskip

\nt and each $\Delta_i$ is irreducible of codimension $\left( {k-2i+2}
\atop 2\right)$ in ${\bf P}^n$.

\bigskip

If one takes
$F_i$ to be an $i\times i$ minor of $\phi$ (in cases (a) and (b)...we'll do 
case (c) later), then $F_i$ has degree $i$ and multiplicity $i-1$
along $Y$.
The strategy we take here for constructing the $F$ to use in Proposition 2.1 is 
therefore to sum general linear combinations of minors of the largest size
until we hit an obstruction (determined
by the corresponding degeneracy locus), 
then to decrease the size of the minor and continue.
Amazingly (at least, to the author), 
we will finish with a log canonical divisor $F$ of multiplicity $r$ 
along $Y$ and degree $e + n$ where $e$ is
independent of $k$ and $m$.

\medskip

We need to invoke two aspects of
the theories of complete linear maps and quadrics (see, for
example, 
\cite{L} for an exposition and specific references).   

\medskip

\nt {\bf Complete Linear Maps and Quadrics...Classical
Construction:} In  both cases (a) and (b), let $U = {\bf P}^n -
\Delta_k$. The space $P$ of  complete objects is the
smooth, projective variety obtained as the  closure of the graph
of $U$ under the morphism:
$$\wedge := (\wedge^2,...,\wedge^k)$$
{\bf Explanation:} A point $\overline \alpha\in {\bf P}^n$ is a
linear map  (modulo scalars). In case (a), it is represented by
a map 
$\alpha: V \rightarrow W^*$, while in case (b), it is represented
by a  symmetric map $\alpha:V \rightarrow V^*$. Thus 
$\wedge^i\alpha$ is a map from $\wedge^iV$ to $\wedge^iW^*
\cong (\wedge^iW)^*$ in (a),
and a symmetric map from $\wedge^iV$ to $\wedge^iV^* \cong (\wedge^iV)^*$ in (b). The 
locus $U \subset {\bf P}^n$ is  the set of $\overline \alpha$ such that  
$\wedge^i\alpha \ne 0$ for all $i \le k$, thus it is where the map:

\medskip

$\wedge  : {\bf P}(V\otimes W) - -> {\bf P}(\wedge^2V 
\otimes \wedge^2W)\times ... \times {\bf P}(\wedge^kV \otimes 
\wedge^kW) \ \mbox{in (a), or}$

\medskip

$\wedge  : {\bf P}(S^2(V)) - -> {\bf P}(S^2(\wedge^2V)) 
\times ... \times {\bf P}(S^2(\wedge^kV)) \ \mbox{in (b)}$

\medskip
  
\nt is defined, and $P$ is embedded in a product of $k$
projective  spaces. 

\bigskip

\nt {\bf Remarks:} $P$ comes equipped with projection morphisms:

\medskip

$\rho_i: P \rightarrow  {\bf P}(\wedge^iV \otimes \wedge^iW)\ \mbox{in (a), and}$
$\rho_i: P \rightarrow {\bf P}(S^2(\wedge^iV))\ \mbox{in case (b).}$

\medskip

\nt Thus, provided $i < k$ in (b) or $i \le k$ and $i < m$ 
in (a), the projection  $\rho_i$
maps to a positive-dimensional projective space, 
giving rise to a base-point-free linear series on $P$. By
definition, the restriction of this linear series to $U$
is spanned by the $i\times i$ minors of $\phi$
(principal minors in case (b)).
To pin down the linear series on $P$ associated to the map
$\rho_i$, we use a second construction of complete linear maps
and quadrics, due to Vainsencher:

\bigskip

\nt {\bf Complete Linear Maps and Quadrics...Blow-Up
Construction:} Recall that we set $X := \ \mbox{bl}({\bf
P}^n,Y)$. The fact  that $i$ matrices of rank one sum to a
matrix of rank at most $i$ implies that in cases (a) and (b),
$$\Delta_i = \Sigma_{i-1}(Y),$$
where $\Sigma_i(Y)$ is the secant variety defined as the closure of 
the union 
of projective planes spanned by $i$ distinct points of $Y$.

\medskip

One blows up 
the degeneracy loci as follows:

\medskip

$f_2: X_2 \rightarrow X_1 := X \ \mbox{blows up the
strict transform of $\Sigma_2(Y) = \Delta_3$}$,

\medskip

$f_3:X_3 \rightarrow X_2 \ \mbox{blows up the strict
transform of $\Sigma_3(Y) = \Delta_4$}$

\medskip

\hskip .5in $\vdots$

$f_{k-1}:\widetilde X := X_{k-1} \rightarrow X_{k-2}\
\mbox{blows up  the strict transform
of $\Sigma_{k-1}(Y)$}$.  

\medskip

\nt For consistency, let $f_1:X \rightarrow {\bf P}^n$ also
be the blow-down. Then:

\medskip

\nt {\bf Theorem:} (\cite{V},Theorem 6.3) (a) 
Each strict transform of $\Delta_{i+2}$
in $X_i$ is smooth, so in particular $\widetilde X$ is
smooth, because:
$$f = f_2 \circ f_3 \circ ... \circ f_{k-1}:
\widetilde X \rightarrow X$$ is a sequence of blow-ups
along smooth centers.  In addition, if we  let $E_i$
denote the strict transform in $\widetilde X$ of the 
(smooth, irreducible) $f_i$-exceptional divisor (for
all $1 \le i \le k-1$), then the divisor:
$$E_1 + ... + E_{k-1}$$
is a normal crossings divisor on $\widetilde X$.

\medskip

(b) The inclusion $\iota: U \hookrightarrow P$ extends to an 
{\bf isomorphism} $\overline \iota:\widetilde X \stackrel \sim\rightarrow P$.

\medskip

\nt (There is also a precise recursive description of $\overline \iota$ which we
will not need here.)

\bigskip

Suppose now that $A_i \subset {\bf P}^n$ is the hypersurface cut out
by some $i \times i$ minor of $\phi$ (or principal minor in case (b)). 
Then $A_i$ has degree $i$ and 
its multiplicity along $\Delta _j$ is (at least) $i - j + 1$ for
all $j \le i$. Thus, $A_i$ determines a divisor:
$$B_i \equiv iH - (i-1)E_1 - ... - E_{i-1}$$
on $P$ by subtracting 
$i-j+1$ copies of $E_j$ from $f^*(A_i)$.

\medskip

I claim that the projection
$\rho_i$ determines a base-point-free sub-linear series 
of $|B_i|$. To see this, 
it suffices to show that the generic multiplicity
of $A_i$ along $\Delta_j$ is precisely $i - j + 1$ 
(so that no $E_j$
is in the base locus of $|B_i|$). But given $A_i$, choose
$\alpha\in \Delta _j - \Delta_{j-1}$ so that 
some $j\times j$ minor (or principal minor) of $\alpha$ contained
in the $i\times i$ minor defining $A_i$ has nonzero 
determinant. Then it is
immediate that $A_i$ has multiplicity exactly $i-j+1$ at
$\alpha$. 

\bigskip

Here, then, is our main proposition to cover cases (a) and (b):

\medskip

\nt {\bf Proposition 3.1:} Given nonnegative integers
$n_2,...,n_k$, let 
$A_{i,1}^g,...,A_{i,n_i}^g$ be the zero loci of general linear
combinations of the determinants of 
$i\times i$ minors of $\phi$ (principal in case (b)).
Let $F_{i,j}^g$ be the strict transform of $A_{i,j}^g$ in $X$,
and let $F^g = \sum_{i=2}^k(F^g_{i,1}+ ... + F^g_{i,n_i})$. Then:

\medskip

(a) In case (a), suppose either $k < m$ or $k = m$ and $n_k
\le 1$.  Then $f:\widetilde X \rightarrow X$ is a log
resolution of 
$(X,F^g)$, and 
the discrepancy at $E_j$ is: 
$$(k-j)(m-j) -1 - \sum_{i > j}(i-j)n_i$$
for each $2 \le j \le k-1$.

\medskip

(b) In case (b), suppose $n_k \le 1$. Then $f:\widetilde X \rightarrow X$
is a log resolution of $(X,F^g)$,
and the discrepancy at $E_j$ is 
$$\left( {k-j+1}\atop 2\right) -1 - \sum_{i > j}(i-j)n_i$$
for each $2 \le j \le k-1$.

\medskip

{\bf Proof:} Vainsencher's theorem tells us the exceptional divisors
have normal crossings. In case (b), and if $k = m$ in case (a),
the last exceptional divisor $E_{k-1}$ is itself the strict
transform of $A_k$. Otherwise,  the strict transforms of the 
$F_{i,j}^g$ in $\widetilde X$ are smooth members
of $|B_i|$, intersecting each other and the 
exceptional divisors transversely by (ordinary) Bertini. 
Thus 
$f$ is a log resolution of the pair $(X,F^g)$.

\medskip

The discrepancies are computed using Riemann-Hurwitz (and the count
for the codimensions at the beginning of this section) as well as the 
linear series computation 
for $|B_i|$, which yields $\sum_{i>j}(i-j)E_j =
f^*(F_{i,j}^g) - f_*^{-1}(F_{i,j}^g)$. 

\bigskip

The following corollary picks out the optimal choices for 
the $n_i$ in order to produce an $e_Y$ which is as small as 
possible. 

\bigskip

\nt {\bf Corollary 3.2:} (a) In case (a), let $n_k = m-k+1$ and 
$n_i = 2$ for $2 \le i < k$. Then $(X,F^g)$ is lc and 
$F^g \equiv (n-1)H - rE$.

\medskip

(b) In case (b), let $n_i = 1$ for all $i$. Then $(X,F^g)$ is lc
and  $F^g \equiv nH - rE$.

\medskip

\nt So using Proposition 2.1, we get $e_Y \le -1$ in case (a) and 
$e_Y \le 0$ in case (b).

\medskip

{\bf Proof:} A direct application of the Proposition tells us that
the discrepancies are all $-1$ for these choices of the $n_i$, so
$(X,F^g)$ is lc. These
and the other computations (the coefficients of $H$ and $E$) are 
straightforward, and left to the reader. 

\bigskip

For case (c), we need versions of the classical construction
and blow-up construction for complete skew
forms.  Specifically, we'll prove the theorem below in an
appendix to this paper:

\medskip

\nt {\bf Complete Skew Forms ... A ``Classical'' Construction:}
Given $V$, let $l = [\frac 12\mbox{dim}(V)]$, let 
$U = {\bf P}(\wedge^2V) - \Delta_{l}$, and consider
the rational map 
$$\wedge: {\bf P}(\wedge^2V) --> {\bf P}(\wedge^4V) \times ...
\times {\bf P}(\wedge^{2l}(V));$$
$$\alpha \mapsto (\alpha\wedge \alpha,\alpha\wedge \alpha\wedge
\alpha,...).$$

It is straightforward to check that for each 
$i=2,...,l$, the degeneracy locus $\Delta_i$ is the 
locus of indeterminacy of the map ${\bf P}(\wedge^2V)
--> {\bf P}(\wedge^{2i}(V))$ obtained by composing 
$\wedge$ with the projection, so in particular,  
$\wedge$ is regular on $U$, and we define:

\medskip

\nt {\bf Definition:} The closure $P := \overline{\Gamma}
\subset {\bf P}(\wedge^2(V)) \times...\times {\bf
P}(\wedge^{2l}(V))$ of the graph of $\wedge$ restricted to 
$U$ is
the space of {\bf complete skew forms} on $V$.

\medskip

\nt {\bf Remark:} This $\wedge$ map is not the same as the 
map we obtain by regarding a skew form as a linear map
and restricting the wedge map.
Firstly, it does not involve the odd wedge powers of $V$,
and secondly it is a ``square root'' of the even part of the 
wedge map in 
the following sense. Notice that $\wedge^{2i}\wedge^2V \subset 
\mbox{Sym}^2\wedge^{2i}V$ as representations of GL$(V)$,
so we can think of $\wedge^{2i}\alpha$ as being a quadratic 
form on $\wedge^{2i}V$. The value of this quadratic form
on a decomposable wedge 
(i.e. a point of the Grassmannian $G(2i,V^*)$)
is always a square, as it is the determinant of a principal 
(skew) minor of the skew form $\alpha:V\rightarrow V^*$ with the 
Pfaffian as a square root. (This is NOT to say, however, that each
$\wedge^{2i}\alpha$ is of rank one.) One checks that these
Pfaffians of 
$2i\times 2i$ principal minors give the linear series associated to 
the rational maps:
$${\bf P}(\wedge^2V) --> {\bf P}(\wedge^{2i}(V)$$
defined above. The common zero scheme of 
these Pfaffians is reduced (unlike the principal determinants), equal to 
the degeneracy locus $\Delta_i$. These linear series will be used 
as before to construct a log canonical divisor. 
  
\bigskip

\nt {\bf Complete Skew Forms ... 
The Construction by Blowing Up:} Let 
$X = \ \mbox{bl}({\bf P}(\wedge^2V),G(2,V^*))$ and blow up the other
degeneracy loci in order as before, letting $X_1 := X$, inductively letting:
$$f_i:X_i\rightarrow X_{i-1} \ \mbox{be the blow up of the 
strict transform of}\ \Delta_{i+1},$$
and letting $f = f_2 \circ ... \circ f_{l-1}:\widetilde X \rightarrow X$. 
Finally, let $E_i \subset \widetilde X$ be the strict transform of 
the $f_i$-exceptional divisor for all $1\le i \le l-1$. Then:

\medskip

\nt {\bf Theorem 3.3:} (a) Each strict transform of
$\Delta_{i+2}$ in 
$X_i$ is smooth, so $\widetilde X$ is smooth, and moreover,
$E_1 + ... + E_{l-1}$ is a normal crossings divisor on $\widetilde X$.

\medskip

(b) The inclusion $\iota:U \hookrightarrow P$ extends to an {\bf isomorphism}
$\overline \iota:\widetilde X \stackrel \sim \rightarrow P$.

\medskip

\nt {\bf Proof:} See the appendix.

\bigskip

Then as before, we conclude that for each $i < \frac 12{\mbox{dim}(V)}$, 
the zero loci $A_i$ of the Pfaffians of principal $2i\times 2i$ minors
of $\phi$ give elements of the linear series:
$$|B_i| := |iH - (i-1)E_1 - ... - E_{i-1}|$$
on $\widetilde X$ which contains  the base-point-free linear
series associated to 
the projection 
$\rho_i:P \rightarrow {\bf P}(\wedge^{2i}(V))$.
Finally, we obtain analogues of 3.1 and 3.2 in case (c):

\medskip

\nt {\bf Proposition 3.4:} Let $n_2,...,n_l$ be 
nonnegative integers, let
$A_{i,1}^g,...,A_{i,n_i}^g$ be zero loci of general linear
combinations of  Pfaffians of $2i\times 2i$ principal minors of
$\phi$,  let
$F_{i,j}^g$ be the  strict transform of $A_{i,j}^g$ in $X$, and
let $F^g = \sum_{i=1}^l(F_{i,1}^g + ... + F_{i,n_i}^g)$. 

\medskip

If $2l < k = \  \mbox{dim}(V)$ or $2l = k$ and $n_l = 1$, then $f:\widetilde X
\rightarrow X$ is a log resolution of $(X,F^g)$ and the discrepancy at $E_j$ is:
$$\left( {k-2j} \atop 2 \right) - 1 - \sum_{i > j}(i-j)n_i$$
for each $2 \le j \le l$.

\medskip

{\bf Proof:} Just as in Proposition 3.1.

\medskip

\nt {\bf Corollary 3.5:} Let $n_l = 1$ if $k$ is even and 
$n_l = 3$ if $k$ is 
odd. Otherwise, let $n_i = 4$ for $2 \le i < l$. Then $(X,F^g)$
is lc and  $F^g \equiv (n-3)H - rE$.

\medskip

\nt Thus using Proposition 2.1, we get $e_Y \le -3$ in case (c).

\newpage

\nt {\bf 4. Curves.} Let $C$ be a smooth, irreducible projective curve
over the complex numbers of genus $g$, let $K_C$ be a canonical
divisor, and
$D$ be a divisor of degree $d \ge 3$. This restriction on $d$
assures us that the linear series map:
$$\phi_{|K_C+D|}: C \rightarrow |K_C+D| \cong {\bf P}^{d+g-2}$$
is an embedding (and we will set $Y = C$ and $n = d+g-2$ in this section).
Notice that the degree of the embedding is $d + 2g-2$, not $d$.

\medskip

We set $X$ to be the blow-up of ${\bf P}^n$ along $Y$ as before,
and recall the standard result (see, for example \cite{ACGH})
that if $d \ge 4$, then the embedded curve $C$ is a scheme-theoretic
intersection of quadric hypersurfaces, from which it follows that:

\medskip

(a) $d_Y = 2$, and

\medskip

(b) the ample cone of $X$ is spanned (in the $H,E$-plane) by $H$ and $2H - E$.

\medskip

The following result is not standard. It is proved in 
\cite{B2} using Thaddeus' stable pairs (\cite{T}) and 
the author's blow-up of secant varieties (\cite{B1}). Indeed, these 
techniques closely resemble the two constructions of the spaces of 
complete objects in the previous section!

\bigskip

\nt {\bf Proposition 4A:} There exist log canonical divisors on $X$ 
that are:

\medskip

(a) numerically equivalent to $(d-1)H - (d-3)E$ (if $g \ge 0$)  

\medskip

(b) numerically equivalent to $dH - (d-2)E$ (if $g > 0$)

\medskip

(c) numerically equivalent to $(\frac{d+g-5}{d-4})(dH - (d-2)E)$
(if $g > 0$ and $d > 4$).

\bigskip

The divisors in (a) and (b) can be taken to be
the strict transforms of hypersurfaces in ${\bf P}^n$,
but in case (c), one obviously needs to stick with ${\bf Q}$-coefficients.

\bigskip

We can apply our strategy directly now in genus $0$ and $1$:

\bigskip

\nt {\bf Proposition 4.1:} If $g = 0$ or $g = 1$, then $e_Y \le 1$.

\medskip

{\bf Proof:} Recall that by Proposition 2.1, we are searching for log 
canonical divisors $F \equiv (n+1)H - (n-1)E$ on $X$ (and here 
$n = d+g-2$). But this 
is just what Proposition 4A (a) and (b) produce for us in genus $0$ and $1$. 

\bigskip

To handle higher genus, we need to improve Proposition 2.1 a bit.

\bigskip 

Namely, thanks to (b) above, we know that as soon as $\epsilon < \frac 12$,  
then  $H - \epsilon E$ is in the the ample cone of 
$X$. Recall that
the key point of Proposition 2.1 was the
observation that the desired vanishing occurs when $p$ and $k$ satisfy:
$$pH - kE \equiv K_X + (1-\epsilon)F + A$$
where $0 < \epsilon < 1$, $F$ is log canonical (at least along $E$) 
and $A$ is big and nef. 
Since Proposition 4A (c) gives us a ``very efficient''
log canonical ${\bf Q}$-divisor on $X$, we'll use this 
divisor and
our better knowledge of the ample cone to get better vanishing results.
 
\bigskip

Assume throughout that $d > 4$
and that  
a log canonical divisor $F$ is given satisfying $F \equiv
(\frac{d+g-5}{d-4})(dH - (d-2)E)$. Then
$$K_X + F \equiv \left(\frac {2g-2}{d-4}\right)(2H - E) - E,$$
which means that if we rewrite $pH - kE$ as above and let 
$\epsilon ' = \frac{d+g-5}{d-4}\epsilon$, then
$$A \equiv (p + d\epsilon ')H - 
(k - 1 + (d-2)\epsilon ')E - \left(\frac {2g-2}{d-4}\right)(2H - E).$$
Thus we see
that $A$ is big and nef if the following {\bf two} inequalities are satisfied,
and at least one of them is strict:

\medskip

(i) $p + d\epsilon ' \ge 2(k - 1 + (d-2)\epsilon ')$ and

\medskip

(ii) $k - 1 + (d-2)\epsilon '\ge \frac {2g-2}{d-4}$.

\bigskip

We now get the following proposition by choosing $\epsilon '$ carefully:

\medskip

\nt {\bf Proposition 4.2:} If $d > \frac {2g+8}3$, then $e_Y \le 1$.

\medskip

{\bf Proof:} For $\epsilon ' \le \frac 1{d-4}$, condition (i) is satisfied
whenever $p \ge 2k - 1$. If there were no additional conditions on 
$k$, then this would imply $e_Y \le 1$. (Recall that $d_Y = 2$).
If we choose $\epsilon '$ to be very close to (and less than) $\frac 1{d-4}$,
then condition (ii) becomes: $k > \frac{2g-d}{d-4} + 1$.

\medskip

By a theorem of Castelnuovo, $C \subset |K_C+D|$ is projectively normal,
which is to say that vanishing holds when $k = 1$ and $p \ge 1$, so that 
we only need to prove vanishing for $k \ge 2$. But moreover, 
by the following:

\medskip

\nt {\bf Proposition (Rathmann):} \cite{R} If $d \ge 5$, then
$$H^i({\bf P}^n,{\cal I}_C^2(p)) = 0\ \mbox{for all}\ i > 0
\ \mbox{and} \ p \ge 3.$$

\medskip

\nt we only need to prove vanishing when $k \ge 3$, which we get since 
both conditions are satisfied if $k \ge 3$ and $d > \frac {2g+8}3$,
hence $A$ is ample, and log Kodaira vanishing applies.

\medskip

\nt {\bf Observation:} If vanishing is proven for 
$p \ge 2k-1$ and $k \le k_0$, then there will be a corresponding 
improvement in the lower bound for $d$ in Proposition 4.2. However, these 
will all be linear in $g$, while I suspect the correct lower bound 
is actually $d \ge 5$, for which I submit the following ``gap'' as evidence.

\medskip

\nt {\bf Proposition 4.3:} If $d \ge 5$ and $g$ are fixed, then 
with at most finitely many exceptions for the values of $p$ and $k$,
$$H^i({\bf P}^{d+g-2},{\cal I}_C^k(p)) = 0
\ \mbox{for all}\ i > 0, k \ge 0 \ \mbox{and} \ p \ge 2k-1.$$

\medskip

{\bf Proof:} By the proof of Proposition 4.2 and the two Propositions cited 
therein, any exceptions must lie in the region in the
$k,p$-plane bounded  on the left by $k = 3$ and on the right by 
$k = \frac {2g-d}{d-4} + 1$. This is of course infinite because
there is no upper bound on $p$.  But notice that if $p \ge 2k -
1 + a$, then we may boost 
$\epsilon '$ to $\frac a{d-4}$ (up to a maximum of $\frac
{d+g-5}{d-4}$), preserving the inequality in condition (i). This
gives  a corresponding lowering of the upper bound for $k$ to
$k < \frac {2g-d - a(d-2)}{d-4} + 1$ beneath which the exceptions
may occur. Thus any exceptions are constrained to lie in a 
roughly triangular region of the plane.

\newpage

\nt {\bf Appendix. Two constructions of complete skew forms.}
Let me begin by arguing why complete skew forms 
(definition in \S 3) are entirely analogous to complete 
quadrics. 

\medskip

Suppose $(R,m)$ is a discrete valuation ring over a field $k$
with quotient field $K$, residue field $k$ and 
uniformizing parameter $t$, and let
$\alpha$ be an $R$-valued skew $2$-form on a vector space $V$
over $k$ of dimension $2l$ or $2l+1$, with the following properties:

\medskip

$\bullet$ the induced ``generic'' $2$-form on $V\otimes_kK$ is of 
rank $2l$, and

\medskip

$\bullet$ the induced ``special'' $2$-form on $V$ is nonzero.

\medskip

In other words, suppose $\alpha$ is the lift of a morphism
$f:\mbox{Spec}(R) \rightarrow {\bf P}(\wedge^2 V)$
with the property that $\wedge\circ f:
\mbox{Spec}(R) \rightarrow {\bf P}(\wedge^2V) \times ...
\times {\bf P}(\wedge^{2l}V)$ is defined 
at the generic point of Spec$(R)$. I want to investigate 
the extension of $\wedge \circ f$ across the special point.

\medskip

Choosing a basis $x_1,...,x_n$ for 
$V^*$ gives a straightforward description of the extension,
since:
$$\alpha = \sum_{i < j} a_{i,j} x_i \wedge x_j$$
with $a_{i,j} \in R$, and an $r$-fold wedge $\alpha \wedge... \wedge \alpha$
is of the form $\sum_I f_I x_{i_1} \wedge...\wedge x_{2r}$
with $f_I \in R$. 
There will be a maximal $d_r$ 
such that each $f_I$ is divisible by $t^{d_r}$, and 
the nonzero images of $t^{-d_r} \alpha \wedge... \wedge \alpha$
in $\wedge^{2r}V^* \otimes R/m \cong \wedge^{2r}V^*$ 
for each $r$
will give the image of the special point.

\medskip

The basis-free approach sets up the analogy with 
complete quadrics. Given $\alpha$, let 
$\alpha_0 \in \wedge^2V^* \otimes R/m\cong \wedge^2V^*$ be its
residue modulo $m$, and suppose $\alpha_0$ has rank $2r_1$. Then
$\alpha_0$ is induced from a nondegenerate skew form on a 
quotient $V \rightarrow T_{r_1}$, and we let 
$W_{r_1}$ be the kernel of this map.
The image of $\alpha$ in $\wedge^2W^*_{r_1}\otimes R$ lies in
$\wedge^2W^*_{r_1}\otimes m$, and we let 
$\alpha_1 \in \wedge^2 W^*_{r_1} \otimes
m/m^2 \cong \wedge^2W^*_{r_1}$ be the residue modulo $m^2$. 
Continuing
in this manner, we produce from $\alpha$ the following data:

\medskip

(D1) A flag of (strict) subspaces:
$$W_{r_m} \subset W_{r_{m-1}} \subset ... \subset W_{r_1}
\subset W_{r_0} = V$$ 
such that dim($W_{r_{i-1}}/W_{r_{i}}) = 2r_{i}$
and $r_1 + ... + r_m = l$.

\medskip

(D2) Skew $2$-forms $\alpha_i$ on $W_{r_i}$, induced from 
nondegenerate skew forms on the quotients $W_{r_i}/W_{r_{i+1}}$.

\medskip

The data (D1) and (D2) determine elements of 
each ${\bf P}(\wedge^{2r}V)$ as follows. For each $\alpha_i$,
let $\widetilde \alpha_i \in \wedge^2V^*$ be an arbitrary lift,
and take the $r$-fold wedge product:
$$\omega_r := \alpha_0 \wedge ... \wedge \alpha_0 \wedge \widetilde 
\alpha_1 \wedge... \wedge \widetilde \alpha_1 \wedge 
\widetilde \alpha_2 \wedge...$$
taking up to $r_1$ copies of $\alpha_0$ followed by up 
to $r_2$ copies of $\widetilde \alpha_1$, etc. until
$r$ terms in all have been taken. It is now an easy exercise
to see that:

\medskip

(i) $\omega_r \in \wedge^{2r}V^*$ is nonzero if $r \le l$, and

\medskip

(ii) $\omega_r$  does not depend upon the choice of lifts, and

\medskip

(iii) if the data (D1) and (D2) come from 
$\alpha \in \wedge^2V^* \otimes R$
as described above, then the $\overline{\omega_r} \in {\bf P}
(\wedge^{2r}V)$ extend  
$\wedge \circ f$ across the special point of Spec$(R)$.

\medskip

In fact, we have the following Lemma, which should look
familiar to anyone who has thought about complete quadrics:

\bigskip

\nt {\bf Lemma A1:} The map to 
${\bf P}(\wedge^2V) \times ...\times {\bf P}(\wedge^{2l}V)$
is a bijection from the set of data (D1) and (D2) 
(modulo scalars)
to the set of complete skew forms.

\medskip

{\bf Proof:} Given a subspace $W \subset V$ of 
codimension $2r$ and quotient $T$, the canonical inclusion
$\wedge^{2r}T^* \otimes \wedge^2W^* \subset \wedge^{2r+2}V^*$
is the key to recovering $\alpha_0,...,\alpha_m$ (modulo scalars) 
from its image in ${\bf P}(\wedge^2V) \times...
\times {\bf P}(\wedge^{2l}V)$. Precisely, if  
$(\beta_1,...,\beta_l)$
is the image, then $\beta_1 = \alpha_0$ (up to scalar multiple)
and therefore determines 
$W_{r_1}$. Then $\beta_{r_1+1}$ determines 
$\alpha_1$ by the inclusion above, 
determining $W_{r_2}$, and $\beta_{r_1+r_2+1}$ determines
$\alpha_2$, etc. proving injectivity.

\medskip

One can always ``smooth'' $\alpha_0,...,\alpha_m$,
taking $\alpha = \alpha_0 + t\widetilde \alpha_1 + ...+
t^m\widetilde\alpha_m$ (where $\widetilde \alpha_i$ denotes
a lift to $\wedge^2V^*$) exhibiting the image of the 
sequence of $\alpha_i$ as a specialization of 
$\wedge \circ f$ and proving that the set of data (D1)
and (D2) maps to the 
complete skew forms. Surjectivity
follows from the valuative criterion for properness.  

\bigskip

The basic idea behind Theorem 3.3 is the same as in the case of
complete quadrics. The normal bundle to the smooth 
subvariety
$\Delta_r - \Delta_{r-1} \subset {\bf P}(\wedge^2V)$ 
is naturally a (twisted) bundle of skew forms on the distinguished
subspaces $W_r \subset V$. Thus, the information consisting
of a rank $2r$ form $\alpha_0$ and a point in the projectivized normal 
bundle to $\Delta_r - \Delta_{r-1}$ at $\overline \alpha_0$
is part of the data (D1) and (D2). Blowing up the degeneracy loci 
in order turns out to give a natural variety structure to 
the data (D1) and (D2) (modulo scalars) which one then proves is 
isomorphic to the variety of complete skew forms.

\medskip

We will thus need to 
consider complete skew forms in 
a relative setting. It will suffice for our purposes
to generalize the above discussion to the case where 
$V$ is a vector bundle over a smooth 
base scheme $X$ over an algebraically closed field, 
though as in the case of  complete bilinear
forms (see \cite{KT}) much of what is proved here can 
presumably be further generalized.

\medskip

Given the vector bundle $V$, 
we introduce the following cast of characters:

\medskip

\nt {\bf Definitions:} (i)
$\pi_{r,V}:{\bf P}(\wedge^{2r}V) \rightarrow X$
 (whenever $2r \le \ \mbox{rk}(V)$).

\medskip

(ii) $\rho_{r,V}: G(V,2r) \rightarrow X$, the 
bundle of $2r$-dimensional {\bf quotients} of the fibers of 
$V$ (with relative Pl\"ucker embedding 
$G(V,2r) \subset {\bf P}(\wedge^{2r}V)$).

\medskip

(iii) $0\rightarrow S_{r,V} \rightarrow 
\rho_{r,V}^*V \rightarrow Q_{r,V}
\rightarrow 0$, the universal sequence on $G(V,2r)$. 

\medskip

(iv) $f_{r,V}:{\bf P}(\wedge^2Q_{r,V}) 
\rightarrow {\bf P}(\wedge^2V)$ inducing a $2$-form on
$V$ from one on the quotient. This is 
the Pl\"ucker embedding when $r = 1$. Note that
${\bf P}(\wedge^2Q_{r,V})$ is a projective bundle over 
$G(V,2r)$, with projection map $\pi_{1,Q_{r,V}}$.

\medskip

(v) $\Delta_{r,V} \subset {\bf P}(\wedge^2V)$ is the ``bundle
of degernacy loci'' of the fibers, well-defined since rank 
is independent of basis. 

\medskip

All the important identifications are made in the following 
lemma.

\medskip

\nt {\bf Lemma A2:} (a) $\Delta_{r,V} =
f_{r,V}({\bf P}(\wedge^2Q_{r,V}))$.

\medskip

(b) $f_{r,V}$ is an embedding when restricted to the 
complement of $\Delta_{r-1,Q_{r,V}}$. The normal bundle of 
the embedding is  $\pi^*\wedge^2 S^*_{r,V}
\otimes {\cal O}(1)$.

\medskip

(c) For each $r\le s$, the fiber product of 
$f_{r,V}$ and $f_{s,V}$ satisfies: 
$${\bf P}(\wedge^2 Q_{r,V}) \times_{{\bf P}(\wedge^2V)}
{\bf P}(\wedge^2 Q_{s,V}) \cong 
{\bf P}(\wedge^2Q_{r,Q_{s,V}})$$ 
which is a projective bundle over the 
flag variety Fl$(V,2s,2r)$. The map to ${\bf P}(\wedge^2Q_{s,V})$
is $f_{r,Q_{s,V}}$ (given by (iv) above) 
and if $\sigma:
\mbox{Fl}(V,2s,2r) \rightarrow
G(V,2r)$ is the forgetful map, then the map
to ${\bf P}(\wedge^2 Q_{r,V})$ is the projection 
from ${\bf P}(\wedge^2Q_{r,Q_{s,V}}) \cong
\sigma^*{\bf P}(\wedge^2Q_{r,V})$.

\medskip

(d) The map from the pull-back of the conormal bundle of 
$f_{r,V}$ to the conormal bundle of $f_{r,Q_{s,V}}$, is, 
after the identifications from (b), the natural map:
$$\sigma^*\left( \wedge^2 S_{r,V}\right) (-1)
\rightarrow
\wedge^2 S_{r,Q_{s,V}}(-1)$$
on the complements of $\Delta_{r-1,*}$ in the fiber product from (c).

\bigskip

{\bf Proof:} A skew form of rank $\le r$ on
a fiber of $V$ is always induced from a skew form on 
an $r$-dimensional quotient of the fiber. This gives (a).

\medskip

In (b), injectivity is clear. Via the Euler sequences for 
the tangent bundles to 
${\bf P}(\wedge^2 V)$ and ${\bf P}(\wedge^2 Q_{r,V})$, one
obtains a sheaf map:
$$0 \rightarrow S_{r,V}^* \otimes Q_{r,V}^* \stackrel \phi\rightarrow
\wedge^2 V^*(1)/\wedge^2Q_{r,V}^*(1)$$
over ${\bf P}(\wedge^2 Q_{r,V})$
with the property that $f_{r,V}$ is an immersion with normal 
bundle isomorphic to the cokernel of $\phi$ wherever $\phi$ is 
fiberwise injective. But $\phi$ factors through the natural map:
$$S_{r,V}^* \otimes Q_{r,V} \rightarrow S_{r,V}^* \otimes Q_{r,V}^*(1)$$
in the obvious way, exhibiting $f_{r,V}$ as an immersion with 
desired normal bundle precisely
on the complement of $\Delta_{r-1,Q_{r,V}}$.

\medskip

(c) is straightforward. The proposed fiber product embeds
naturally in the product. And (d) follows from a diagram chase.

\bigskip

Now that the identifications (a)-(d) in Lemma A2 have been
established, it is a formal consequence of the recursive nature
of the conormal bundles and maps among them that:

\medskip

(1) The blow up, in order, of the strict transforms of the 
degeneracy loci $\Delta_{r,V}$ produces a smooth variety
$\widetilde X$ with normal crossings exceptional divisors, and

\medskip

(2) The set of data (D1) and (D2) (modulo scalars) corresponds
to the points of $\widetilde X$, with each new subspace
$W_{r_i}$ and skew form $\alpha_i$ corresponding to 
a normal direction of the strict transform
of $\Delta_{r_i,V}$ (modulo scalars).

\bigskip

This is proved, for instance, in Proposition 2.2 of \cite{B1} with 
Lemma 1.3 of that paper playing the role of Lemma A2 here. The 
steps of that proof are readily adapted to handle this
situation (or, for that matter, the case of complete linear 
maps and quadrics). This takes care of the first part of Theorem
3.3, leaving us to  prove that the rational map:
$$\iota: {\bf P}(\wedge^2 V) --> {\bf P}(\wedge^2 V)
\times ... \times {\bf P}(\wedge^{2l}V)$$
extends to an embedding of $\widetilde X$ which agrees with 
the map on the set of data (D1) and (D2) (modulo scalars)
considered in Lemma A1. 

\medskip

To prove that $\iota$ extends to a map from $\widetilde X$,
we use:

\medskip

\nt {\bf Lemma A3:} Suppose $X$ is a normal variety
and $Y$ is a projective variety over an algebraically 
closed field $k$. If $f:X -->Y$ is a rational map
and $\overline f:X \rightarrow Y$
extends $f$ as a map of sets, then $\overline f$
is a morphism if and only if the following ``valuative''
criterion is satisfied:

\medskip

(*) For all discrete valuation rings $R$ over $k$
with residue field $k$, and all morphisms
$\alpha:\ \mbox{Spec}(R) \rightarrow X$ sending
the generic point $\xi \in \ \mbox{Spec}(R)$ 
to the domain of $f$,
the image of the special point under 
$\overline f \circ \alpha$
agrees with the specialization of $f(\alpha(\xi))$.

\medskip

{\bf Proof:} This is an immediate consequence of 
Zariski's Main Theorem applied to the graph of $f$
(see \cite{H},V.5.2). 

\bigskip

We can apply this to the extension of $\iota$ 
via the identification of 
$\widetilde X$ with the set of data (D1) and (D2)
(modulo scalars). The discussion preceding Lemma A1 
tells us that the conditions of Lemma A3 are 
satisfied, implying that the bijection of 
Lemma A1 is a morphism $\overline \iota: \widetilde
X \rightarrow P$ to the space of complete
skew forms. Thus it only remains to prove that 
$\overline \iota$ is an immersion. But we can prove
this by induction on the rank of $V$ (again 
considering the relative setting). Namely, we
know that $\overline \iota$ is an immersion on 
the complement of exceptional divisors, since that
locus is included in ${\bf P}(\wedge^2V)$. On the other
hand, by induction, the exceptional 
divisor over $\Delta_{r,V}$ embeds in 
the complete skew forms
on ${\bf P}(\wedge^2 Q_r)$, hence in 
the complete skew forms on $V$ via the embedding:
$${\bf P}(\wedge^2 Q_r) \times_{G(V,2r)}
... \times {\bf P}(\wedge^{2r-2}Q_r)
\hookrightarrow {\bf P}(\wedge^2 V) \times_X 
 ... \times {\bf P}(\wedge^{2r-2}V) \times G(V,2r).$$ 

This only leaves normal vectors to the exceptional
divisros to worry about. But such a normal vector 
is either tangent to some other exceptional 
divisor, and we have already dealt with it, or else it
maps to a nonzero normal vector to the smooth 
part of a $\Delta_{r,V} \subset {\bf P}(\wedge^2 V)$
under the blow-down. Thus in all cases, 
nonzero tangent vectors to 
$\widetilde X$ remain nonzero under $\overline \iota_*$,
and $\overline \iota$
is indeed an embedding.

University of Utah, Salt Lake City, UT 84112

{\it email address:} bertram@math.utah.edu

\end{document}